\documentclass[12pt]{article}
\newtheorem{ttt}{Theorem}
\newtheorem{ddd}{Definition}[section]

\newtheorem{lll}{Lemma}[section]

\newcommand{\aut}{automorphism\hspace{2pt}}

\title{The group of automorphisms of semigroup End$(P[X])$}
\author{A. Berzins\\University of Latvia\\ e-mail: aberzins@latnet.lv}
\date{August 2004}

\begin{document}
\maketitle

AMS \textit{Mathematics Subject Classification}: 14A99, 14P05


\begin{abstract}

In this paper is proved that the group of automorphisms of
semigroup End$(P[X])$, if $P$ is algebraically closed field, is
generated by semi-inner automorphisms.

\end{abstract}

\section{Introduction}

The basis of the classical algebraic geometry is the Galois
correspondence between \textit{P}-closed ideals in
$P[x_1,x_2,\ldots,x_n]$ and algebraic sets in the affine space
$P^n$. Having noted that a point $a=(a_1,a_2,\ldots,a_n)$ unique
determines a homomorphism $x\in\textrm{Hom}(P[X])$
($s(x_i)=a_i$), B. Plotkin in [7-10] makes foundation of the
universal algebraic geometry, i.e., algebraic geometry for
arbitrary variety $\Theta$ of universal algebras. Let us fix a
variety $\Theta$ and an algebra $H\in\Theta$. Let $W=W(X)$ be a
free algebra in $\Theta$ with finite set of generators. The set
$\textrm{Hom}(W,H)$ we consider as affine space of points over
$H$. Then arbitrary congruence $T$ in $W$ determines the set of
points in affine space $\textrm{Hom}(W,H)$:
$$T'_H=A=\{\mu\in\textrm{Hom}|T\subset\textrm{Ker}\mu\}.$$
A set of points $A\subset\textrm{Hom}(W,H)$ determines congruence
of W
$$A'=A'_W=T=\bigcap_{\mu\in A}\textrm{Ker}\mu.$$
We call the set of points $A$ such that $A=T'$ for some $T$ an
algebraic set in $\textrm{Hom}(W,H)$. A relation $T$ with $T=A'$
for some $T$ is a congruence in $W$. We call such a congruence an
$H$-closed one.

It was found [3-4, 7-10] that many problems in the universal
algebraic geometry, such as geometric equivalence, geometric
similarity, isomorphism and equivalence of categories of
algebraic sets and varieties depend on structure of ${\rm
Aut}(\Theta^0)$ and Aut(End($W$)), where $\Theta^0$ is the
category of algebras $W=W(X)$ with the finite $X$ that are free
in $\Theta$, and $W$ is a free algebra in $\Theta$. Structure of
${\rm Aut(Com}-P))^0$ (i.e, the classical case) was described in
[2]. Here was introduced variant of the concept of quasiinner
automorphism - basic  concept for description of Aut(End($W$)) and
$\textrm{Aut}(\Theta^0)$. Later groups $\textrm{Aut}(\Theta^0)$
were described for categories of free Lie algebras and free
associative algebras [1, 5, 6]. Some other definition of
quasiinner automorphism was used in mentioned papers. Also are
described ${\rm Aut}(\Theta^0)$ for varieties of groups,
semigroups and some another varieties. For some varieties is
described Aut(End($W$)), but these problems are more difficult.
In particular, this problem was not solved even in the classical
case of free commutative algebras, i.e. the structure of
Aut(End$(P[x_1,\ldots ,x_n])$ is not describe if $n>2$. In all
solved cases groups Aut(End($W$)) are generated by semi-inner
automorphisms and mirror automorphism. Now in presented paper the
problem in classical case is solved.

I want to note that the proof in the present paper does not use
the description of the group ${\rm Aut}(P[X])$. It is very
important, because structure of the group ${\rm Aut}(P[x_1,\ldots
,x_n])$ is not described, when $n>2$, but counterexamples show
that structures of ${\rm Aut}(P[x_1,\ldots ,x_n])$ and ${\rm
Aut}(P[x,y]$) are principally different. So, methods presented in
the paper may be useful to describe the group Aut(End($W_n$)) for
different varieties of algebras over the field.

\section{Definitions}

Let us recall some definitions for variety of commutative algebras
$\textrm{Com}-P$ (see \textrm{[4, 10]} for general case).

\begin{ddd}

Let $P[X])=P[x_1,\ldots ,x_n]$ be a free commutative algebra over
a field $P$ with a finite set of generators and
$\tau\in\textrm{Aut}(\textrm{End}P[X])$. It is known
$\textrm{[4]}$ that there exists a bijection $\mu:P[X]\mapsto
P[X])$ such that for every $s\in \textrm{End}(P[X])$
$$ s^{\tau}=\mu s\mu^{-1}.$$
Such representation of $\tau$ is called a representation of $\tau$
as quasiinner automorphism generated by $\mu$.

\end{ddd}
Of course, arbitrary set of substitutions does not generate an
automorphism, But there are two kinds of bijections $\mu$ that
generate an automorphism.\\

1. $\mu=\overline{\alpha}$, where $\overline{\alpha}$ is the
natural extension of an automorphism $\alpha \in
\textrm{Aut}(P[X])$ on $W$, i.e. $\overline{\alpha}\left( \sum
a_ku_k\right)=\sum \alpha(a_k)u_k$.
 We shall write $\overline{\alpha}=\alpha$ and call it "\aut of
 field".\\

 2. $\mu=\eta\in \textrm{Aut}(P[X])$.\\

\begin{ddd}

The automorphism $\tau$ generated by an automorphism
$\mu\in\textrm{Aut}(P[X])$ is called inner automorphism.

\end{ddd}

\begin{ddd}

Note, that every $\alpha \in \textrm{Aut}(P[X])$ belongs to
normalizer of subgroup  $\textrm{Aut}(P[X])$ in the group of all
bijections of $P[X])$. So, every product of elements of
$\textrm{Aut}P[X])$ and $\textrm{Aut}(P([X])$ can be represented
in the form $\mu=\alpha\eta$, $\alpha \in \textrm{Aut}P[X])$,
$\eta\in\textrm{Aut}(P[X])$. Automorphism
$\tau\in\textrm{Aut}(\textrm{End}(P[X])$ generated by $\mu$ is
called semiinner.

\end{ddd}

Let $W(X)$ be a free finitely generated algebra in arbitrary
variety $\Theta$.

\begin{ddd}

Bijection $\mu : W(X)\mapsto W(X)$ is called central, if for every
$s\in \textrm{End}(W(X))$
$$\mu s=s\mu.$$
Algebra $W(X)$ is called central, if every its central bijection
is identical.

\end{ddd}

\begin{ttt}

Let $W(X)$ be a free central algebra and $\mu$ some bijection of
$W(X)$ generating an automorphism
$\tau\in\textrm{Aut}(\textrm{End}(W(X)))$. Then $\mu$ transforms
every base of algebra $W(X)$ to a base of algebra.

\end{ttt}
\textbf{Proof.} [see 4]\\[6pt]

For free algebras in a variety $\Theta$ of algebras over field $P$
(Com-P, Ass-P, Lie-P and others) we shall correct the definition
of central algebra. Clearly, in every such algebra $W(X)$ we have
a bijection $\mu(u)=au+b$, where $a,b\in P,$ $a\neq 0$ (linear
bijection). This bijection commutate with every endomorphism of
algebra $W(X)$ and, obviously, transform every base of $W(X)$ to
a base.

\begin{ddd}

Finitely generated algebra $W(X)$ over field $P$ is called almost
central, if every its central bijection is linear.

\end{ddd}

\begin{ttt}

Free finitely generated commutative, associative are almost
central.

\end{ttt}
\textbf{Proof.} Let $\mu$ be a central bijection of $W(X)$,
$\mu(x_i)=r_i(x_1,x_2,\ldots,x_n)\in W(X)$ , and $s\in
\textrm{End}(W(X)$, $s(x_1)=x_1$, $s(x_i)=0$ for $i=1$. Then
$$\mu s(x_1)=\mu (x_1)=r_1(x_1,x_2,\ldots,x_n),$$
$$s\mu(x_1)=s(r_1(x_1,x_2,\ldots,x_n))=r_1(x_1,0,\ldots,0)=r(x_1).$$
So $r_1(x_1,x_2,\ldots,x_n)=r(x_1)$ is a polynomial of one
variable.

For arbitrary $u\in W(X)$ take $s\in \textrm{End}(W(X))$
$s(x_1)=u$, $s(x_i)=x_i$ for $i>1$. Then
$$\mu(u)=\mu s(x_1)=s\mu(x_1)=s(r(x_1))=r(u),$$
and so $\mu(u)=r(u)=a_0+a_1u+\cdots+a_ku^k$. Because $\mu$ is
bijection, so $r$ is linear, Q.E.D..\\[3pt]
\textbf{Corollary.} \textit{From the theorems 1 and 2 we have that
in mentioned algebras every bijection $\mu$, which generate an
automorphism of $\textrm{End}(W(X))$  transform every base of
algebra $W(X)$ to a base of algebra.}

\section{The group ${\rm Aut(End(P[X]))}$}

Let $\tau$ be an automorphism of $\textrm{End}(P[X])$ generated by
bijection $\mu$. We call $s\in \textrm{End}(P[X])$ a constant
endomorphism if $\textrm{Im}(s)=P$. Clearly, if $s$ is constant
then $s^{\tau}$ also is constant. Denote the set of all constant
endomorphisms by $Const$.

Let $a\in P$ and $s$ be a constant endomorphism such that
$s(x_1)=a$. Then $\mu(a)=\mu(s(x_1))=s^{\tau}(\mu(x_1))\in P$. So
restriction of $\mu$ on $P$ is a bijection $s:P\mapsto P$.

Let $\mu(0)=a$ and $\mu(1)=b$. We shall consider the linear
bijection $l:P[X])\mapsto P([X])$, $l(u)=cu+d$ such that $l(a)=0$
and $l(b)=1$. Since linear bijection generates identical
automorphism of the semigroup $\textrm{End}(P[X])$, so bijections
$\mu$ and $\mu'=l\mu$ generate equal automorphisms. Note that
$\mu'(0)=0$ and $\mu'(1)=1$. We denote $\mu=\mu'$ and assume that
$\mu(0)=0$ and $\mu(1)=1$.

Since $\mu$ transforms base of algebra $P[X]$ to base, so the
endomorphism $\eta:$ $\eta(x_i)=\mu(x_i)$ is an automorphism of
$P[X]$. Now we shall consider the bijection $\mu'=\eta^{-1}\mu$.
We have for $\mu'$ that $\mu'(x_i)=x_i$. Let as denote $\mu=\mu'$
and prove that $\mu$ is multiplicative function on $P[X]$.

\begin{lll}
Let $\mu$ be a bijection of $P[X]$ generating an automorphism
$\tau\in \textrm{Aut}(\textrm{End}(P[X]))$ such that $\mu(0)=0$,
$\mu(1)=1$ and $\mu(x_i)=x_i$. Then for every $u,v\in P[X]$ such
equality holds:
$$\mu(u\cdot v)=\mu(u)\cdot\mu(v)$$
\end{lll}
\textbf{Proof.} (This proof use algebraic closure of $P$. It will
be useful give an algebraic proof of this lemma for arbitrary
infinite field)

Note that bijection $\mu$ generating an automorphism determines a
bijection of algebraic sets in $P^n$ of codimension 1. It follows
from Hilbert Nullstellensatz that
$$\mu(x_1x_2)=c\mu(x_1)^k\mu(x_2)^l.$$
We shall write $\mu(x_1x_2)\approx \mu(x_1)^k\mu(x_2)^l$.
Substitution (endomorphism) $s:x_1\rightarrow u, x\rightarrow v$
gives similarity $\mu(uv)\approx\mu(u)^k\mu(v)^l$. Using
substitutions $u\rightarrow u, v\rightarrow1$ we get that
$\mu(u)\approx\mu(u)^k$, so $k=1$. Analogously, $l=1$ and
$\mu(uv)=cx\mu(u)\mu(v)$. Substitution $u\rightarrow 1$,
$v\rightarrow 1$ gives that $c=1$, Q.E.D..

\begin{lll}

With the preceding notation
$$\mu(u+v)=\mu(u)+\mu(v)$$

\end{lll}
\textbf{Proof.} Let $a\in P$, $\mu(x_1)=x_1$,
$\mu(x_1-a)=r(x_1,\ldots,x_n)$. Since $\langle x_1,(x_1-a)\rangle$
is the unit ideal, so $\langle x_1,r(x_1,\ldots,x_n)\rangle$ also
is unit, i.e. the system of equations
$\{x_1=0,r(x_1,\ldots,x_n)=0\}$ is inconsistent. Therefore,
$r(0,x_2,\ldots,x_n)=0$ and $r=d+x_1\cdot g(x_1,\ldots,x_n)$.
Since $r$ is basic element (i.e., the element of some base), so
$r-d$ also is basic. Clearly, basic element is irreducible. So
$x_1\cdot g(x_1,\ldots,x_n)$ is irreducible, i.e.
$g(x_1,\ldots,x_n)=c$ and $r=cx_1+d$.

Let now $a\in P$
\begin{center}
$\mu(x_1-a)=c\mu(x_1)+d$,\\
$\mu(x_1+a)=c_1\mu(x_1)+d_1$.
\end{center}
Thus for every $u\in P[X]$
\begin{center}
$\mu(u-a)=c\mu(u)+d$,\\
$\mu(u+a)=c_1\mu(u)+d_1$.
\end{center}
So $0=\mu(a-a) =c\mu(a)+d$, $d=-c\mu(a)$.\\
From equality $1=\mu(1)=\mu((-1)\cdot(-1))=\mu(-1)^2$ we have
$\mu(-1)=-1$ and $\mu(-a)=-\mu(a)$. So
$0=\mu(-a+a)=-c_1\mu(a)+d_1$. Thus
\begin{center}
$\mu(u-a)=c(\mu(u)-\mu(a))$,\\
$\mu(u+a)=c_1(\mu(u)+\mu(a)),$
\end{center}
$\mu(u)=\mu(u-a+a)=c_1(\mu(u-a)+\mu(a))=c_1c\mu(u)+(-c_1c+c_1)\mu(a)$.
So $cc_1=1,c_1(c-1)=0$ and $c=c_1=1$. Thus in the field $P$ we
have an identity
$$\mu(a+b)=\mu(a)+\mu(b).$$
It means that restriction of $\mu$ on $P$ is an automorphism of
$P$. Let $\mu'=\epsilon^{-1}\mu$. Then $\mu'$ is identical on
$P$, and
\begin{equation}
\mu'(au)=a\mu'(u)
\end{equation}
for every $a\in P$ and $u\in P[X]$. Considering all constant
endomorphisms we get
$$\left(\mu(x_1+x_2)-\mu(x_1)- \mu(x_2)\right)\in\bigcap_{s\in Const}\textrm{Ker}s=0.$$
Substitution $s(x_1)=u,s(x_2=v$ gives required identity.\\
Now we can formulate the main result of this paper.

\begin{ttt}
The group $\textrm{Aut}(\textrm{End}(P[X]))$, where $P$ is
algebraically closed field, is generated by semiinner
automorphisms.
\end{ttt}
It follows from lemma 3.1, lemma 3.2 and identity (1).\\[6pt]

\noindent{\Large\textbf{Acknowledgments}}\\[12pt]
The author is happy to thank professor B. Plotkin for stimulating
discussions of the results.\\[30pt]

\noindent{\Large\textbf{References}}
\begin{enumerate}

\item A. Berzins, {\it The group of automorphisms of the category
of free associative algebras} (to appear)

\item A. Berzins, {\it The automorphisms of \,{\rm End}{}$K[x]$}, Proc.
Latvian Acad. Sci., Section B, 2003, vol. 57, no. 3/4, pp. 78-81.

\item A. Berzins, {\it Geometric equivalence of algebras}, Int. J.
Alg. Comput., 2001, vol. 11, no. 4, pp. 447-456.

 \item A. Berzins, B. Plotkin and E. Plotkin, {\it Algebraic
geometry in Varieties with the Given Algebra of constants}, J.
Math. Sci., New York, 2000, vol. 102, no. 3, pp. 4039-4070.

\item R. Lipyanski, B. Plotkin {\it Automorphisms of categories
 of free modules and free Lie algebras}, Preprint

\item G. Mashevitzky, B. Plotkin, E. Plotkin, {\it Automorphisms
of categories of free Lie algebras.}, Journal of Algebra, Vol.
282, 2004, pp. 490-512.

\item B. Plotkin, {\it Algebraic logic, varieties of algebras and
algebraic varieties}, Proc. Int. Alg. Conf., St. Petersburg, 1995,
Walter de Gruyter, New York, London, 1996.

\item B. Plotkin, {\it Varieties of algebras and algebraic
varieties}, Israel J. of Mathematics, 1996, vol. 96, no. 2, pp.
511-522.

\item B. Plotkin, {\it Varieties of algebras and algebraic
varieties. Categories of algebraic varieties}, Sib. Adv. Math.,
1997, vol. 7, no. 2, pp. 64-97.

\item B. Plotkin, {\it Algebras with the same (algebraic)
geometry}, Proc. Steklov Inst. Math, Vol. 242, 2003, pp. 165-196.

\end{enumerate}

\end{document}